\documentclass[a4paper,11pt]{amsart}

\usepackage[T1]{fontenc}
\usepackage{amssymb,amsfonts,amsmath,amsthm}
\usepackage{url}
\usepackage{color}
\usepackage{setspace}
\usepackage{enumerate}
\usepackage{epsfig,graphicx,psfrag}
\usepackage{asymptote,tikz}
\usepackage{pinlabel}

\input{xy}
\xyoption{all}
\objectmargin={3mm}

\newcounter{notes}
\newtheorem{theorem}{Theorem}
\newtheorem{proposition}[theorem]{Proposition}

\theoremstyle{definition}

\newtheoremstyle{theoremwithref}{}{}{\itshape}{}{\bfseries}{.}{.5em}{#1 #2 #3}
\theoremstyle{theoremwithref}

\newcommand{\C}{\mathbb{C}}
\newcommand{\R}{\mathbb{R}}
\newcommand{\Z}{\mathbb{Z}}

\newcommand{\SL}{\mathrm{SL}}
\newcommand{\PSL}{\mathrm{PSL}(2,\mathbb{R})}

\newcommand{\Hom}{\mathrm{Hom}}
\newcommand{\HH}{\mathbb{H}^2}
\newcommand{\Isom}{\mathrm{Isom}}
\newcommand{\tr}{\mathrm{tr}}
\newcommand{\Homeop}{\mathrm{Homeo}_+(S^1)}
\newcommand{\rot}{\mathrm{rot}}

\title{Dynamical rigidity of non discrete representations in $\PSL$}

\author{Maxime Wolff}
\address{Sorbonne Universit\'es, UPMC Univ. Paris 06, Institut de Math\'ematiques
de Jussieu-Paris Rive Gauche, UMR 7586, CNRS, Univ. Paris Diderot,
Sorbonne Paris Cit\'e, 75005 Paris, France}
\email{maxime.wolff@imj-prg.fr}

\begin{document}

\begin{abstract}
The aim of this note is to advertise on a result, not stated explicitely,
but proved, in \cite{wol08}. Namely, if $\Gamma$ is any group, if $\rho_1,\rho_2$
are representations of $\Gamma$ in $\PSL$, one of them being non elementary
and non discrete, and if for all $\gamma\in\Gamma$, $\rho_1(\gamma)$
and $\rho_2(\gamma)$ have the same rotation number, then $\rho_1$ and
$\rho_2$ are conjugate in $\PSL$. In particular, if two non discrete,
non elementary representations yield semi-conjugate actions on the circle,
then they are conjugate in $\PSL$.
\end{abstract}

\maketitle


Recall that the group $\PSL$ of orientation-preserving isometries of the
hyperbolic plane $\HH$ may be considered as a subgroup of the group
$\Homeop$ of orientation-preserving homeomorphisms of the circle
(the circle at infinity of~$\HH$). Any element of
$\Homeop$ has a rotation number, invariant under semi-conjugation
(see eg \cite{ghy01,Katie} for definitions, and a lot of information).
If $A\in\PSL$, its rotation number, $\rot(A)\in\R/2\pi\Z$ is equal to zero
if $A$ is hyperbolic, parabolic or the identity, and to its angle of
rotation if $A$ is elliptic.

If $\Gamma$ is a fundamental group of a compact orientable surface,
any Teichm\"uller representation $\rho$ of $\Gamma$ in $\PSL$ admits
continuous deformations (in the Teichm\"uller space),
which are $C^0$-conjugate, but not even
$C^1$-conjugate to $\rho$ (see \cite{ghy85}, Proposition III.4.1).
More generally, if $\Gamma$ is any finitely generated group, then every non-elementary,
discrete representation $\rho$ of $\Gamma$ in $\PSL$ may be deformed
to representations $C^0$-conjugate to $\rho$, but not conjugate in $\PSL$,
by deforming the Fuchsian group $\rho(\Gamma)$.
In contrast, we may observe the following.

\begin{proposition}\label{LaProp}
Let $\Gamma$ be any group, and $\rho_1,\rho_2\in\Hom(\Gamma,\PSL)$.
Suppose $\rho_1$ is non elementary, non discrete, and suppose that
for every $\gamma\in\Gamma$, $\rot(\rho_1(\gamma))=\rot(\rho_2(\gamma))$.
Then there exists $g\in\PSL$ such that $\rho_1=g\rho_2 g^{-1}$.
\end{proposition}

In the statement above, a representation $\rho$ is called elementary if
$\rho(\Gamma)$ is elementary, ie, has a finite orbit in $\HH\cup\partial_\infty\HH$
(see eg \cite{Katok}, section 2.4).

Although Proposition~\ref{LaProp} was proved in~\cite{wol08}
(see Theorem 2.39 in \cite{wol08}), it was not stated explicitely there,
nor mentionned in the introduction. It was even not included in the
published version, \cite{wol11}.
The aim of this note, thus, is to make this statement more accessible.
It is motivated, in part, by the recent article
\cite{KKM},
as Proposition~\ref{LaProp} answers a part of
Question~1.13 of \cite{KKM}.

Even though the proof of Proposition~\ref{LaProp} is already in \cite{wol08},
it is short and elementary so I reproduce it here.

\begin{proof}
By Selberg's lemma, the subgroup $\rho_1(\Gamma)$ of $\PSL$ admits a
torsion-free finite index subgroup, which, thus, is still non elementary
and non discrete. Hence, $\rho_1(\Gamma)$ admits an elliptic element of
infinite order (see \cite{Katok}, Theorem~2.4.5), say, $\rho_1(\gamma_0)$, for
some $\gamma_0\in\Gamma$.
We may conjugate $\rho_1$ so that
$\rho_1(\gamma_0)=\pm\left(\begin{array}{cc}\cos\theta & \sin\theta \\
-\sin\theta & \cos\theta\end{array}\right)$,
with $\theta\in(0,\pi)$
irrational (indeed, this matrix is a rotation of angle $2\theta$ in $\HH$).
Also, $\rho_2(\gamma_0)$, having rotation number $2\theta$, is elliptic, conjugate to $\rho_1(\gamma_0)$.
Up to conjugating $\rho_2$, we may further suppose that $\rho_1(\gamma_0)=\rho_2(\gamma_0)$.

Now let $\gamma\in\Gamma$. For $i\in\{1,2\}$, write
$\rho_i(\gamma)=\pm\left(\begin{array}{cc}a_i & b_i \\ c_i & d_i\end{array}\right)$.
Then for all $n\in\Z$,
\[ \left|\tr\rho_i(\gamma\gamma_0^n)\right|=\left|(a_i+d_i)\cos(n\theta)+(c_i-b_i)\sin(n\theta)\right|;
\]
this is the absolute value of the scalar product of the vectors
$(a_i+d_i,c_i-b_i)$ and $(\cos(n\theta),\sin(n\theta))$ of $\R^2$.
The first vector cannot be zero, and
the second vector may be chosen as close as we want from any point of the unit circle.
Thus, for infinitely many values of $n$, the quantity
$|\tr(\rho_1(\gamma\gamma_0^n))|$ is in $(0,2)$, hence $\rho_1(\gamma\gamma_0^n)$
is elliptic, thus, so is $\rho_2(\gamma\gamma_0^n)$, and the cosines of their
rotation numbers are equal, hence
\[|(a_1+d_1)\cos(n\theta)+(c_1-b_1)\sin(n\theta)|=|(a_2+d_2)\cos(n\theta)+(c_2-b_2)\sin(n\theta)|.
\]
This equality, valid for infinitely many distinct values of $n$, implies that
$|a_1+d_1|=|a_2+d_2|$. This proves that $|\tr(\rho_1(\gamma))|=|\tr(\rho_2(\gamma))|$,
for all $\gamma\in\Gamma$.

It is well known that, if the traces of two
irreducible representations in $\SL(2,\C)$ agree on each element of the group,
then they are conjugate, see \cite{CS}, Proposition 1.5.2.
An adaptation of this argument (see eg \cite{wol11}, Proposition 2.15)
gives, in our context, that $\rho_1$ and $\rho_2$ are conjugated by
an element $g$ of $\Isom(\HH)$. This $g$ preserves the orientation of the
hyperbolic plane,
otherwise we would have $\rot(\rho_1(\gamma_0))=-\rot(\rho_2(\gamma_0))$.
Finally, $\rho_2=g\rho_1 g^{-1}$, with $g\in\PSL$.
\end{proof}

\small

Acknowledgements. Proposition~\ref{LaProp} arose from a question that my
advisor Louis Funar asked me shortly after my PhD defence.

\normalsize


\end{document}